\documentclass[11pt, oneside]{amsart}
\usepackage{graphicx}
\usepackage{amssymb}
\usepackage{amsfonts}
\usepackage{amscd}
\usepackage{amstext}
\usepackage{epsfig}
\usepackage{mathrsfs}
\usepackage{amsmath,amsthm,amsfonts,amssymb,amscd}

\newtheorem{teo}{Theorem}[section]
\newtheorem{cor}[teo]{Corollary}
\newtheorem{lem}[teo]{Lemma}
\newtheorem{pro}[teo]{Proposition}
\newtheorem{defi}[teo]{Definition}

\newtheorem{prob}{Problem}[section]

\newcommand{\N}{{\mathbb{N}}}
\newcommand{\Z}{{\mathbb{Z}}}

\newcommand{\R}{{\mathbb{R}}}
\newcommand{\dom}{{\mbox{dom}\,}}

\newcommand{\bbox}{\ \rule[-0.3mm]{2.5mm}{2.5mm}}

\newcommand{\no}{\noindent}
\newcommand{\xx}{{\mbox{\hspace{.2cm}{\LARGE -\kern -.47em $x$}}}}

\newcommand{\mm}{\medskip}
\newcommand{\bb}{\bigskip}

\begin{document}

\title[$C^r-$closing for flows on surfaces]{On $C^r-$closing for flows on orientable and non-orientable 2--manifolds}

\author{Carlos Gutierrez and Benito Pires}

\address{Instituto de Ci\^encias Matem\'aticas e de Computa\c{c}\~{a}o, Universidade de
S\~ao Paulo, Caixa Postal 668, 13560-970, S\~ao Carlos--SP, Brazil.}

\email{gutp@icmc.usp.br, \,\,\, bfpires@hotmail.com}

\maketitle

\begin{abstract} We provide an affirmative answer to the $C^r-$Closing Lemma, $r\ge 2$,
for a large class of flows defined on every closed surface.
\end{abstract}

%%%%%%%%%%%%%%%%%%%%%%%%%%%%%%%%%%%%%%%%%%%%%%%%%%%%%%%%%%%%%%%%%%%%%%%%%%%%%%%%%%%%%%%%%%%%%%%%%%%%%%%%%%%%%%%%%
\section{Introduction}

This paper addresses the open problem $C^r$ Closing Lemma, which can
be stated as follows:

\begin{prob}[$C^r$ Closing Lemma]\label{prob1} Let $M$ be a compact smooth manifold, $r\ge 2$ be an integer,
$X\in\mathfrak{X}^r(M)$ be a $C^r$ vector field on $M$, and $p\in M$
be a non--wandering point of $X$. Does there exist
$Y\in\mathfrak{X}^r(M)$ arbitrarily $C^r-$close to $X$ having a periodic trajectory passing through $p$?
\end{prob}

C. Pugh \cite{Pg1} proved the $C^1$ Closing Lemma for flows and
diffeomorphisms on manifolds. As for greater smoothness $r\ge 2$,
the $C^r$ Closing Lemma is an open problem even for flows on the
$2-$torus. Concerning flows on closed surfaces, only a few, partial
results are known in the orientable case (see \cite{car,G4,G11}). No
affirmative $C^r-$closing results are known for flows on
non--orientable surfaces. In this paper, we present a class of flows
defined on every closed surface supporting non--trivial recurrence
for which Problem \ref{prob1} has an affirmative answer -- see
Theorem A. Notice that every closed surface distinct from the
sphere, from the projective plane and from the Klein bottle \mbox{(see
\cite{MK})} admits flows with non--trivial recurrent trajectories
(see \cite{Hi}).

To achieve our results we provide a partial, positive answer to the
following local version of the $C^r$ Closing Lemma for flows on
surfaces:

\begin{prob}[Localized $C^r$ Closing Lemma]\label{prob2} Let $M$ be a closed surface, $r\ge 2$ be an integer,
$X\in\mathfrak{X}^r(M)$ be a $C^r$ vector field on $M$, and $p\in M$
be a non--wandering point of $X$. For each neighborhood $V$ of $p$ in $M$
and for each neighborhood $\mathcal{V}$ of $X$ in
$\mathfrak{X}^r(M)$, does there exist $Y\in\mathcal{V}$, with $Y-X$ supported in $V$, having a periodic
trajectory meeting $V$?
\end{prob}

It is obvious that if Problem \ref{prob2} has a positive answer for
some class of vector fields $\mathscr{N}\subset\mathfrak{X}^r(M)$
then so does Problem \ref{prob1}, considering the same class
$\mathscr{N}$. The approach we use to show that a flow has local
$C^r-$closing properties is to make arbitrarily small
$C^r-$twist--perturbations of the original flow along a transversal
segment. This requires a tight control of the perturbation: it may
happen that a twist--perturbation leaves the non--wandering set
unchanged \cite{GB1} or else collapses it into the set of
singularities \cite{car}, \cite{G6}. More precisely: \mbox{C.
Gutierrez} \cite{G6} proved that local $C^2-$closing is not always
possible even for flows on the \mbox{$2$--torus;} C. Carroll
\cite{car} presented a flow on the $2$--torus with poor
$C^r-$closing properties: no arbitrarily small
$C^2-$twist--perturbation yields closing; \mbox{C. Gutierrez and B.
Pires \cite{GB1}} provided a flow on a non--orientable surface of
genus four whose non--trivial recurrent behaviour persists under a
class of arbitrarily small $C^r-$twist--perturbations of the
original flow.

Deeply related to Problem \ref{prob1} is the Peixoto--Wallwork
Conjecture that the Morse-Smale vector fields are $C^r-$dense on
non--orientable closed surfaces, which is implied by the following
open problem:

\begin{prob}[Weak $C^r$ Connecting Lemma]\label{prob3}
Let $M$ be a non--orientable closed surface, $r\ge 2$ be an integer,
and $X\in\mathfrak{X}^r(M)$ have singularities, all of which
hyperbolic. Assume that $X$ has a non--trivial recurrent trajectory.
Does there exist $Y\in\mathfrak{X}^r(M)$ arbitrarily $C^r-$close to
$X$ having one more saddle--connection than $X$?
\end{prob}

M. Peixoto \cite{Pe} gave an affirmative answer to the Weak $C^r$
Connecting Lemma, $r\ge 1$, for flows on orientable closed surfaces
whereas C. Pugh \cite{Pu67} solved the Peixoto--Wallwork Conjecture in
class $C^1$.

%The only
%obstruction to extending Peixoto's proof for non--orientable compact
%surfaces is the existence of flows with ``non--orientable"
%non--trivial recurrences. The existence of such flows was first
%shown by C. Gutierrez \cite{G10} on the torus with two cross--caps
%(a non--orientable surface of genus four). Afterwards \mbox{A.
%Nogueira \cite{No2}} showed that these kinds of flows abound on any
%non--orientable closed surface of genus greater or equal to four.

%Besides the $C^r-$Closing Lemma type results that we have already
%mentioned, we would like to recall some other relevant results in
%this topic. Namely, we should mention Ma\~n\'e's $C^1$--Ergodic
%Closing Lemma~\cite{M2}, Herman's $C^r$--Closing
%Lemma~\cite{He1,He2}, Hayashi's $C^1$--Connecting
%Lemma~\cite{Hay1,Hay2,Hay3}. In the same way as above, we wish to
%mention \cite{B-C,Ol1,Ol2,Pg2,Pg3,Pg4,P-R}.

\medskip

To give a positive answer to the Peixoto--Wallwork Conjecture, it
would be enough to prove either the $C^r-$Closing Lemma or the Weak
$C^r$ Connecting Lemma for the class $\mathscr{G}^\infty(M)$ of
smooth vector fields having nontrivial recurrent trajectories and
finitely many singularities, all hyperbolic. However there is not a
useful classification of vector fields of $\mathscr{G}^\infty(M).$
Surprisingly, this is not contradictory with the fact that the class
$\mathscr{F}^\infty(M)$ of smooth vector fields having nontrivial
recurrent trajectories and finitely many singularities, all locally
topologically equivalent to hyperbolic ones, is essentially
classified. The vector fields that are constructed to classify
$\mathscr{F}^\infty(M)$ have flat singularities \cite{G2}. The
answer to either of the following questions is unknown (see
\cite{MSMM} for related results):
\medskip

\noindent (1) Given $X\in\mathscr{F}^\infty(M)$, is there a vector
field $Y\in \mathscr{G}^\infty (M)$ topologically equivalent
\mbox{to $X$?}\vspace{0.2cm}

\noindent(2) Given $X\in \mathscr{G}^\infty(M)$ which is dissipative
at its saddles, is there $Y\in\mathscr{G}^\infty(M)$ topologically
equivalent to $X$ but which has positive divergence at some of its
saddles?

\medskip

Considering  vector fields of $\mathscr{G}^\infty(M)$ which are
dissipative at their saddles, their existence in a broad context
was considered by C. Gutierrez \cite{G8}. The motivation of this
work was to find a \mbox{$C^r-$ Closing Lemma} for all vector fields of
$\mathscr{G}^\infty(M)$ whose existence is ensured by the work done
in \cite{G8}. In this paper we have accomplished this aim. We do not
know any other existence result improving that of \cite{G8}.

%%%%%%%%%%%%%%%%%%%%%%%%%%%%%%%%%%%%%%%%%%%%%%%%%%%%%%%%%%%%%%%%%%%%5

\section{Statement of the results}\label{results}

Throughout this paper, we shall denote by $M$ a closed Riemannian
surface, that is, a compact, connected, \mbox{boundaryless,}
$C^\infty$, Riemannian \mbox{2--manifold} and by
$\mathfrak{X}_H^r(M)$ the open subspace of $\mathfrak{X}^r(M)$
formed by the \mbox{$C^r$ vector} fields on $M$ having singularities
(at least one), all of which hyperbolic. When $M$ is neither the
torus nor the Klein bottle, $\mathfrak{X}_H^r(M)$ is also dense in
$\mathfrak{X}^r(M).$ To each \mbox{$X\in\mathfrak{X}_H^r(M)$} we shall
associate its flow $\{X_t\}_{t\in\R}$. Given a transversal segment
$\Sigma$ to $X\in\mathfrak{X}_ H^r(M)$ and an arc length
parametrization $\theta:I\subset\R\to \Sigma$ of $\Sigma$, we shall
perform the identification $\Sigma=\theta(I)=I$, where $I$ is a
subinterval of $\R$. In this way, subintervals of $I$  will denote
subsegments of $\Sigma.$ If $P:\Sigma\to\Sigma$ is the forward
Poincar\'e Map induced by $X$ on $\Sigma$ and $x$ belongs to the
domain ${\rm dom}\,(P)$ of $P$, we shall denote:
$$DP(x)= D(\theta^{-1}\circ P\circ \theta)(\theta^{-1}(x)).$$ Notice
that $DP(x)$ does not depend on the particular arc length
parametrization $\theta$ \mbox{of $\Sigma$} and may take positive
and negative values. Given $n\in\N\setminus\{0\}$, we let
$$\mathcal{O}_n^-(\partial\Sigma)=\{P^{-i}(\partial\Sigma):0\le i\le
n-1\},$$ where $\partial\Sigma$ denotes the set of endpoints of
$\Sigma$ and $P^0$ is the identity map. In this way, the $n-$th
iterate $P^n$ is differentiable on ${\rm
dom}\,(P^n)\setminus\mathcal{O}_n^-(\partial\Sigma)$.

\begin{defi}[Infinitesimal contraction]\label{map} Let $\Sigma$ be a
transversal segment to a vector field \linebreak{$X\in\mathfrak{X}_H^r(M)$} and
let $P:\Sigma\to\Sigma$ be the forward Poincar\'e Map induced by
$X$. Given  \linebreak {$n\in\N\setminus\{0\}$} and {$0<\kappa<1$}, we say
that $P^n$ is an {\it infinitesimal $\kappa$-contraction} if $\vert
D{P^n}(x)\vert<\kappa$ for all \mbox{$x\in{\rm
dom}\,(P^n)\setminus\mathcal{O}_n^-(\partial\Sigma)$}.
\end{defi}

We
say that $N\subset M$ is a {\it quasiminimal set} if $N$ is the
topological closure of a non--trivial recurrent trajectory of $X$.

\begin{defi}\label{defi2} We say that $X \in \mathfrak{X}_H^r (M)$
has the infinitesimal contraction property at a subset $V$ of $M$ if
for every non--trivial recurrent point $p\in V$, for every
$\kappa\in(0,1)$ and for every transversal segment $\Sigma_1$ to $X$
passing through $p$, there exists a subsegment $\Sigma$ of
$\Sigma_1$ passing through $p$ such that the forward Poincar\'e Map
$P:\Sigma\to\Sigma$ induced by $X$ is an infinitesimal
$\kappa$--contraction.
\end{defi}

Given a transversal segment $\Sigma$ to $X\in\mathfrak{X}_H^r(M)$ passing through a non--trivial recurrent point
of $X$,
we let $\mathscr{M}_P(\Sigma)$ denote the set of Borel probability measures on $\Sigma$ invariant by the forward
Poincar\'e Map $P:\Sigma\to\Sigma$ induced by $X$. We say that a Borel subset $B\subset \Sigma$ is {\it of
total measure} if $\nu(B)=1$ for all $\nu\in\mathscr{M}_P(\Sigma)$.

\begin{defi}[Lyapunov exponents] We say that
$X\in\mathfrak{X}_H^r(M)$ has negative Lyapunov exponents at a subset $V$ of $M$ if for
each non--trivial recurrent point $p\in V$ and for each transversal
segment $\Sigma_1$ passing through $p$, there exist a subsegment
$\Sigma$ of $\Sigma_1$ containing $p$ and a total measure set
$W\subset\Omega_+$ such that for all $x\in W$,
\[\chi(x)=\liminf_{n\to\infty}\frac1n \log\vert DP^n(x)\vert<0,
\]
where $P:\Sigma\to\Sigma$ is the forward Poincar\'e Map induced by $X$
and $\Omega_+=\cap_{n=1}^\infty\, {\rm
dom}\,(P^n)$.
\end{defi}

Now we state our results.\newline

\noindent{\rm\bf Theorem A.} {\it Suppose that
$X\in\mathfrak{X}_H^r(M)$, $r\ge 2$, has the contraction property at
a {quasimi\-nimal} set $N$. For each $p\in N$, there exists
\mbox{$Y\in\mathfrak{X}_H^r(M)$} arbitrarily $C^r-$close to $X$ having a periodic trajectory passing through $p$.}
\vspace{0.2cm}

\noindent{\bf Theorem B. }{\it Suppose that $X$ has divergence less or equal to zero at
its saddle--points and that $X$ has negative Lyapunov exponents at a
quasiminimal set $N$. Then $X$ has the infinitesimal contraction
property at $N$.}\vspace{0.2cm}

\noindent{\bf Theorem C.\,} {\it  Suppose that
$X\in\mathfrak{X}_H^r(M)$, $r\ge 2$, has the contraction property at
a {quasimi\-nimal} set $N$.  There exists $Y\in\mathfrak{X}^r(M)$
arbitrarily $C^r-$close to $X$ having one more saddle--connection
than $X$.}\vspace{0.2cm}

%%%%%%%%%%%%%%%%%%%%%%%%%%%%%%%%%%%%%%%%%%%%%%%%%%%%%%%%%%%%%%%%%%%%%%%%%%%%%%%%%%%%%%
\section{Preliminares}

A transversal segment $\Sigma$ to $X\in\mathfrak{X}_H^r(M)$ {\it
passes through $p\in M$} if $p\in\Sigma\setminus\partial\Sigma.$
%and $\Sigma$ is {\it compact} if $\Sigma\supset\partial\Sigma$.
Given $p\in M$, we shall denote by $\gamma_p$ the trajectory of $X$
that contains $p$. We may write $\gamma_p=\gamma_p^-\cup\gamma_p^+$
as the union of its negative and positive semitrajectories,
respectively. We shall denote by $\alpha(p)$ or $\alpha(\gamma_p)$
(resp. $\omega(p)$ or $\omega(\gamma_p)$) the $\alpha-$limit set
(resp. $\omega-$limit set) of $\gamma_p$. The trajectory $\gamma_p$
is {\it recurrent} if it is either {\it $\alpha-$recurrent} (i.e.
$\gamma_p\subset\alpha(\gamma_p)$) or {\it $\omega-$recurrent}
 (i.e. $\gamma_p\subset\omega(\gamma_p)$). A recurrent
trajectory is either {\it trivial} (a singularity or a periodic
trajectory) or {\it non--trivial}. A point $p\in M$ is {\it
recurrent} (resp. {\it non--trivial recurrent,
$\omega-$recurrent,...}) according to whether $\gamma_p$ is
recurrent (resp. non--trivial recurrent, $\omega-$recurrent...). We
say that $N\subset M$ is a {\it quasiminimal set} if $N$ is the
topological closure of a non--trivial recurrent trajectory of $X$.
There are only finitely many quasiminimal sets $\{N_j\}_{j=1}^m$,
all of which are invariant. Furthermore, every non--trivial
recurrent trajectory is a dense subset of exactly one quasiminimal
set.

\begin{pro}\label{pro3} Let $N$ be a quasiminimal set of
$X\in \mathfrak{X}_H^r(M)$. Suppose that for some non-trivial
recurrent point $p \in N$, there exist a transversal segment
$\Sigma$ to $X$ passing through $p$, $(\kappa,n)\in(0,1)\times \N$,
and $L > 0$ such that the forward Poincar\'e Map $P:\Sigma\to\Sigma$
induced by $X$ has the following properties:
\begin{enumerate}
\item[$(a)$] The $n$-th iterate $P^n$ is an infinitesimal $\kappa-$contraction;
\item[$(b)$] $\sup \,\{\vert DP(x)\vert:x\in{\rm dom}\,(P)\} \leq L$.
\end{enumerate}
Then $X$ has the infinitesimal contraction property at $N$.
\end{pro}
\proof We claim that
\begin{itemize}
\item[(a)] for every $K\in(0,1)$  there exists a subsegment $\Sigma_K$
of $\Sigma$ passing through $p$ such that the forward Poincar\'e Map
$P_K:\Sigma_K\to\Sigma_K$ induced by $X$ is an infinitesimal
$K$--contraction.
\end{itemize}
In fact,  let $L_0 = \max\, \{1,L^{n - 1}\}$ and $d \in \N$ be such
that $L_0 \kappa^d < K$. We shall proceed considering only the case
in which  $p$ is nontrivial $\alpha-$recurrent. We can take a
subsegment $\Sigma_K$ of $\Sigma$ passing through $p$ such that $\;
\; \Sigma_K \subset \dom(P^{- dn})\; $ and $\;\Sigma_K,
P^{-1}(\Sigma_K),\cdots, P^{- dn}(\Sigma_K)\;$ are paiwise disjoint.
Hence, if $P_K:\Sigma_K\to\Sigma_K$ is the forward Poincar\'e Map
induced by $X,$ then, for all $q\in\dom(P_K)$, there exists
$m(q)>dn$ such that $P_K(q)=P^{m(q)}(q).$ In this way, since the
function $m:q\mapsto m(q)$ is locally constant, $|DP_K(q)| =
|DP^{m(q)}(q)| \le L_0 \kappa^d< K$ for all $q\in{\rm
dom}\,(P_K)\setminus P_K^{-1}(\partial\Sigma_K)$. This proves (a).

\mm

Let $q\in N$ be a nontrivial recurrent point.
 Now we shall shift the property obtained in (a) to any
segment $\widetilde \Sigma$ transversal to $X$ passing through $q.$
We shall only consider the case in which $q$ is non--trivial
$\alpha-$recurrent and so $\gamma_q^-$ is dense in $N$.

Let $K\in (0,1)$ and take $p_1\in (\gamma_q^-\cap
\Sigma_{K/2})\setminus\{p\}.$ \;
 Select a subsegment $\;\Sigma_1\;$ of $\;\Sigma_{K/2}\;$
passing through $p_1$ and a subsegment $\;\widetilde \Sigma_K\;$ of
 $\;\widetilde \Sigma\;$  passing through $q$ such that the forward Poincar\'e Map $T: \Sigma_1
\to \widetilde \Sigma_K$ is a diffeomorphism and, for all $x\in
\Sigma_1$, $y\in \widetilde\Sigma_K$, \mbox{$|DT(x)DT^{-1}(y)|<2.$}
This implies that the forward Poincar\'e Map $\widetilde P_K:
\widetilde \Sigma_K\to \widetilde \Sigma_K$ will be an infinitesimal
$K-$contraction because

$$|D\widetilde P_K(y)| = | D(T\circ P_1\circ T^{-1}) (y)| \le 2
|DP_1(z)|< K,$$

\no where $P_1:\Sigma_1\to \Sigma_1$ is the forward Poincar\'e Map
induced by $X$ and $T(z)=y.$
\endproof

\begin{defi}[flow box] Let $X\in\mathfrak{X}_H^r(M)$ and let
$\Sigma_1,\Sigma_2$ be disjoint, compact transversal segments to $X$
such that the forward Poincar\'e Map $T:\Sigma_1\to\Sigma_2$ induced
by $X$ is a \mbox{diffeomorphism.} For each $p\in\Sigma_1$, let
$\tau(p)=\min\,\{t>0:X_t(p)\in\Sigma_2\}$.
 The compact region
\mbox{$\{X_t(p):p\in\Sigma_1,\,0\le t\le \tau(p)\}$} is called a
flow box of $X$.
\end{defi}

\begin{teo}[flow box theorem]\label{fbt}
 Let $U\subset M$ be an open set, $X\in\mathfrak{X}_H^r(U)$, \linebreak
 $\Sigma\subset U$ be a compact transversal segment to $X$
and $p \in\Sigma\setminus\partial\Sigma$. There exist $\epsilon>0$ \linebreak arbitrarily small
such that \mbox{$B=B(\Sigma,\epsilon)=\{X_t(p):t\in [-\epsilon,0]\,,
p\in\Sigma\}$} is a flow box of $X$, and a $C^r-$diffeomor\-phism
{$h:B\to [-\epsilon,0]\times [a,b]$} such that $h(p)=(0,0)$,
$h(\Sigma)=\{0\}\times [a,b]$, $h\vert_\Sigma$ is an isometry and
$h_*(X\vert_B)=(1,0)\vert_{[-\epsilon,0]\times [a,b]}$, where
$a<0<b$, $(1,0)$ is the unit horizontal vector field in $\R^2$ and
$h_*(X\vert_B)$ is the pushforward of the vector field $X\vert_B$ by
$h$. The map $h$ is denominated a $C^r-$rectifying diffeomorphism
\mbox{for $B$.}
\end{teo}
\proof See Palis and de Melo \cite[Tubular Flow Theorem, p. 40]{PaMe}. \endproof

\begin{defi} Given a compact transversal segment
$\Sigma$ to \mbox{$X\in\mathfrak{X}_H^r(M)$,}
$p\in\Sigma\setminus\partial\Sigma$ and $\epsilon>0$ small, we say
that $B(\Sigma,\epsilon)=\{X_t(p):t\in [-\epsilon,0]\,,p\in\Sigma\}$
is a  flow box of $X$ ending at $\Sigma$ or at $p$. We say that $B(\Sigma,\epsilon)$ is arbitrarily thin if
$\epsilon$ can be taken arbitrarily small and we say that $B(\Sigma,\epsilon)$ is arbitrarily small if $B(\Sigma,\epsilon)$
can be taken contained in any neighborhood \mbox{of $p$.}
\end{defi}

Next lemma will be used in the proofs of Theorem \ref{saddle-connection} and Theorem \ref{mt1}.

\begin{lem}\label{pro5} Suppose that $X\in\mathfrak{X}_H^r(M)$ has the
infinitesimal contraction property at a non--trivial recurrent point
$p\in M$ of $X$. There exist an arbitrarily small flow box $B_0$ of $X$
ending at $p$ and an arbitrarily small neighborhood $\mathcal{V}_0$
of $X$ in $\mathfrak{X}_H^r(M)$ such that every $Z\in\mathcal{V}_0$,
with $Z-X$  supported in $B_0$, has the infinitesimal contraction
property at $B_0$.
\end{lem}
\proof Let $\Sigma_1=(a_1,b_1)$ be a transversal segment to $X$ passing through $p$
such that the forward Poincar\'e Map $P_1:\Sigma_1\to\Sigma_1$ induced by $X$ is an
infinitesimal $\kappa-$contraction for some $\kappa\in (0,1)$. Let $[a,b]\subset (a_1,b_1)$
be a compact subsegment passing through $p$ and let $B_0=B([a,b],\epsilon)$ be a flow box.
There exists a neighborhood $\mathcal{V}_1$ of $X$ in $\mathfrak{X}_H^r(M)$ such that for every $Z\in\mathcal{V}_1$
with $Z-X$ supported in $B_0$ we have that $B_0$ is still a flow box of $Z$. In particular, for every $Z\in\mathcal{V}_1$
such that $Z-X$ supported in $B_0$, ${\rm dom}\,(P_Z)={\rm dom}\,(P_1)$, where $P_Z$ denotes the forward
Poincar\'e Map induced by $Z$ on $(a_1,b_1)$. Given $\delta>0$ satisfying $0<\kappa+\delta<1$, by the continuity of
the map $Z\mapsto DP_Z$, there exists a neighborhod
$\mathcal{V}_0\subset\mathcal{V}_1$ of $X$ such that for every $Z\in\mathcal{V}_0$ with
$Z-X$ supported in $B_0$ we have that $\vert DP_Z(w)\vert<\vert DP_1(w) \vert+\delta<\kappa+\delta<1$ for all
$w\in {\rm dom}\,(P_1)$. Hence $P_Z$ is an infinitesimal $(\kappa+\delta)$--contraction. The rest of the proof
follows as in the proof of \mbox{Proposition \ref{pro3}} by recalling that the trajectory of every non--trivial recurrent point of $Z$ in $B_0$ meets $(a_1,b_1)$.
\endproof
%Given $0<\overline{\kappa}<\kappa<1$ and
%$\delta=\kappa-\overline{\kappa}$, let $\Sigma\subset V$ be a
%transversal segment passing through $p$ such that $P=P_\Sigma$ is an
%infinitesimal $\overline{\kappa}-$contraction. Let $\epsilon>0$ be
%such that $B=B(\Sigma,\epsilon)$ is a flow box ending at $p$
%contained in $V$. It is plain that if $Z\in\mathcal{N}(B,X)$ is such
%that $Z-X$ is supported in $B$ then ${\rm dom}\,(P_Z)={\rm
%dom}\,(P)$, where $P_Z$ is the forward Poincar\'e Map induced by $Z$
%on $\Sigma$. Since $P_Z$ depends smoothly on $Z$, we have that there
%exists a neighborhood \mbox{$\mathcal{V}_0\subset\mathcal{N}(B,X)$}
%of $X$ in $\mathfrak{X}_H^r(M)$ such that if $Z\in\mathcal{V}_0$ is
%such that $Z-X$ is supported in $B$ then \mbox{$0<\vert D
%{P_Z}(x)\vert<\vert D{P}(x)\vert+\delta <k$}, for every $x\in{\rm
%dom}\,( {P_Z})$.. Hence, ${P_Z}$ is an infinitesimal $k-$contraction
%for each $Z\in\mathcal{V}_0$ such that $Z-X$ is supported in $B$. By
%proceeding as in the proof of \mbox{Proposition \ref{glob}}, we may
%show every $Z\in\mathcal{V}_0$, with $Z-X$ supported in $B$, has the
%infinitesimal contraction property at $B$.

\section{Topological Dynamics}

Let $X\in\mathfrak{X}_H^r(M)$. We say that $N\subset M$ is an {\it
invariant set of $X$} if $X_t(N)\subset N$ for all $t\in \R$. We say
that $K\subset N$ is a {\it minimal set of $X$} if $K$ is compact,
non--empty and invariant, and there does not exist any proper subset
of $K$ with these properties. We shall need the following lemmas
from topological dynamics.

As every vector field of $\mathfrak{X}_H^r(M)$ has singularities,
the Denjoy--Schwartz Theorem \linebreak (see \cite{Sch} or \cite[pp.
39--40]{V}) implies that
\begin{lem}\label{td1} Let $X\in\mathfrak{X}_H^r(M)$, $r\ge 2.$
Then any minimal set of $X$ is either a singularity or a periodic
trajectory.
\end{lem}

The proof of the following lemma can be found in  \cite[Theorem
2.6.1]{NiZh}.

\begin{lem}\label{td2} Let $X\in\mathfrak{X}_H^r(M)$ and let $p\in M$.
Then $\omega(p)$ $($resp. $\alpha(p))$ is exactly one of the
following sets: a singularity, a periodic trajectory, an attracting
graph, or a quasiminimal set.
\end{lem}

\begin{lem}\label{td3} Let $N$ be a quasiminimal set of $X\in\mathfrak{X}_H^r(M)$. Then  every trajectory of $N$ is either a
saddle--point or a saddle--connection or else a non--trivial
recurrent trajectory dense in $N$ $(${\it which may possibly be a
saddle--separatrix.}$)$
\end{lem}
\proof See  \cite[Theorem 2.4.2, pp. 31--32]{NiZh}. \endproof

%\begin{lem}\label{Peixoto1} Let $X \in \mathfrak{X}_H^r(M), r \geq 1$,
% and let $\gamma$ be a nontrivial recurrent trajectory of $X.$
%Suppose that every separatrix contained in $\omega(\gamma)$ is a
%saddle--connection and that $\omega(\gamma)$ contains a
%saddle--point. Then $\omega(\gamma)$ is an attracting graph.
%\end{lem}
%\proof See  \cite{NiZh}. ??? \endproof

%\begin{lem}\label{Peixoto2} Let $X \in \mathfrak{X}_H^r(M), r \geq 2$,
% and let $N$ be a quasiminimal
%set of $X$. Then there exist one or two saddle--points $s_1, s_2 \in
%N$, two saddle--separatrices $\sigma_1, \sigma_2\subset N$, with
%$\omega(\sigma_1)=s_1$ and $\alpha(\sigma_2)=s_2$ such that
%$\alpha(\sigma_1)=N=\omega(\sigma_2)$.
%\end{lem}
%\proof Firstly let us proof that $X$ has singularities in $N$ and
%all of them are hyperbolic saddle--points. If this was not the case,
%then $N$ would contain no singularities and, by Lemma \ref{td3}, $N$
%would be a minimal set of $X$ contradicting Lemma \ref{td1}.

\begin{lem}\label{Peixoto2} Let $X \in \mathfrak{X}_H^r(M), r \geq
2$, and let $N$ be a quasiminimal set of $X$. Then there exist
saddle--separatrices $\sigma_1, \sigma_2 \subset N$ such that
$\alpha(\sigma_1)=N=\omega(\sigma_2)$.
\end{lem}

\proof Firstly let us proof that $X$ has singularities in $N$ and that
all of them are hyperbolic saddle--points. If this was not the case,
then $N$ would contain no singularities and, by \mbox{Lemma \ref{td3},} $N$
would be a minimal set of $X$ contradicting Lemma \ref{td1}.

We shall only prove that $N$ contains dense unstable separatrices.
Suppose by contradiction that

\begin{itemize}
  \item [(a)] every unstable separatrix $\sigma \subset N$
  is a saddle--connection.
\end{itemize}

Take a non-trivial $\omega$-recurrent semitrajectory $\gamma^+$ in $N$ (there is a continuum of such
trajectories in $N$, see \cite[Theorem 2.1, p. 57]{ArBeZh}). We
 say that a region $R\subset M$ is a $\gamma^+$-flow-box  if
there exists a homeomorphism $\;\;h:[-1,1] \times [0,1] \to
R \;\;$ such that

\begin{itemize}
  \item [(b1)] for all $y \in (0,1]$, \;  $h ([-1,1] \times \{y\})$ is
  an arc of trajectory of $X$ starting at the point $h((-1,y))$
  and ending at the point $h((1,y))$. Also, $h((0,0))$ is a
  saddle--point and
 $h([-1,0) \times \{0\})$ (resp. $h((0,1] \times \{0\})$ is a
  stable (resp. an unstable) half--separatrix of $h((0,0))$;
  \item [(b2)] $h(\{-1\} \times [0,1])$ (resp. $h(\{1\} \times
  [0,1])$) is a transversal segment to $X$ called the
  entering edge (resp. exiting edge) of $R$. Moreover, $\gamma^+
  \cap h(\{-1\} \times [0,1])$ accumulates at  the point
  $h((-1,0)).$
\end{itemize}

As $X$ has only finitely many unstable separatrices, by using (a)
we shall be able to find a sequence $\;\; R_1,R_2,\ldots,R_n \;\;$ of
$\gamma^+$-flow-boxes, whose interiors are pairwise disjoint, such
that, for all $i \in \{1,2,\ldots,n-1\}$, the exiting edge of $R_i$
is the equal to the entering edge of $R_{i+1}$ and the exiting edge
of $R_n$ is contained in the entering edge of $R_1$. In this way,
the interior of $\;\;\cup^{n}_{i=1} R_i \;\;$ is an open annulus
eventually trapping the semitrajectory $\gamma^+$ which so cannot be dense. This contradiction proves the lemma.
\endproof

\begin{defi}\label{defi6}
Let $X\in\mathfrak{X}_H^r(M)$ and let $\sigma$ be a non--trivial
recurrent unstable separatrix of a saddle--point $s.$ We say that a
transversal segment $\Sigma$ to $X$ is $\sigma$-adapted if $\sigma$
$($oriented as starting at $s)$ intersects $\Sigma$ infinitely many
times and the first two of such intersections are precisely the
endpoints of $\Sigma.$
\end{defi}

%If $\Sigma$ is a transversal segment  to a  vector field $X$ and
%$\theta:[a,b] \to \Sigma$ is an arc length parametrization of
%$\,\Sigma$, we shall perform the identification
%$\Sigma=\theta([a,b]) = [a,b]$. In this way, subintervals of $[a,b]$
%(whether open or closed) will denote subsegments of $\Sigma$. In the
%proof of Theorem~\ref{saddle-connection}, we shall use the fact that
%$[a,b]$ may be represented by another segment of the form
%$[a+s,b+s].$

\begin{lem}\label{lem7}
Let $\sigma$ be a non--trivial recurrent unstable saddle--separatrix
of \mbox{$X\in\mathfrak{X}_H^r(M)$}. Then every transversal
segment $\Sigma_1=(a_1,b_1)$ to $X$ intersecting $\sigma$ contains a
compact subsegment $[a,b]\subset (a_1,b_1)$ which is
$\sigma-$adapted.
\end{lem}

%\begin{lem}\label{lem7}
%Let $p$ be a nontrivial backward recurrent point of
%$X\in\mathfrak{X}_H^r(M)$ which is accumulated by an unstable
%separatrix $\sigma$ of a saddle--point. Then any open interval
%$(a_1,b_1)$ containing $p$ and transversal to $X$ contains a closed
%subinterval $[a,b]$ intersecting $\gamma_p^-$  which is
%$\sigma-$adapted.
%\end{lem}

\proof Orient $\sigma$ so that it starts at the saddle--point
$\alpha(\sigma).$ Let $p_1, p_2, p_3$ be the first three points at
which $\sigma$ intersects $(a_1,b_1)$ and denoted in such a way that
\mbox{$a_1<p_1<p_2<p_3<b_1.$} If $\sigma$ accumulates at $p_2$ from
above (resp. from below) then $[p_2,p_3]$ (resp. $[p_1,p_2]$) will
be $\sigma-$adapted.
\endproof

\begin{lem}\label{lem8} Let $X\in\mathfrak{X}_H^r(M)$,
$\Sigma=[a,b]$ be a transversal segment to $X$ passing through a non--trivial recurrent point of $X$  and
\mbox{$P:[a,b]\to [a,b]$} be the forward Poincar\'e Map induced by
$X$. Then ${\rm dom}\, (P)\setminus \{a,b\}$ is properly
contained in $(a,b)$ and consists of finitely many open intervals
such that if $s \notin \{a,b\}$ is an endpoint of one of these
intervals then the positive semitrajectory $\gamma^+_{s}$ starting
at $s$ goes directly to a saddle--point without returning to
$[a,b].$
\end{lem}
\proof The proof of this lemma can be found in Palis and de Melo \cite[pp. 144--146]{PaMe} or in Peixoto \cite{Pe}.
\endproof

\section{$C^r-$Connecting Results}

\begin{defi}\label{N} Given $X\in\mathfrak{X}_H^r(M)$ and a flow box $B$ of
$X$, we shall denote by $\mathcal{A}(B,X)$ the set of the vector fields $Y\in\mathfrak{X}_H^r(M)$
supported in $B$ such that for all $\lambda\in [0,1]$, $B$ is still a flow box of $X+\lambda Y$.
\end{defi}

%\begin{lem}\label{Zhuzhoma}
%Let $X\in\mathfrak{X}_H^r(M)$. If a regular point $p\in M$ of $X$ is
%contained in the $\omega-$limit set of a trajectory of $X$ then $p$
%belongs either to a closed orbit or to a saddle connection or else
%to a nontrivial recurrent trajectory.
%\end{lem}
%\no \emph{\textbf{Proof.}} See

%\mm

In next lemma we assume that the domain of the forward Poincar\'e
Map $P$ is non--empty. In the applications of Lemma \ref{lem9} and
Theorem \ref{tt}, $p$ will be a non--trivial recurrent point.

\begin{lem}\label{lem9}
Let $X\in\mathfrak{X}_H^r(M)$ be smooth in a neighborhood $V_0$ of a
point \mbox{$p\in M$} and let \mbox{$\Sigma=[a,b]\subset V_0,$} with
$a<0<b$, be a transversal segment to $X$ passing through $p=0$.
There exist an arbitrarily thin flow box
\mbox{$B=B([a,b],\epsilon)$} contained in $V_0$, and
\linebreak{$Y\in\mathcal{A}(B,X)\subset\mathfrak{X}_H^r(M)$}
arbitrarily $C^r-$close to the zero--vector--field such that for
each \mbox{$\lambda\in [0,1]$} the forward Poincar\'e Map
\mbox{$P_\lambda: [a,b]\to [a,b]$} induced by $X+\lambda Y$ is of the form $P_\lambda=E_\lambda\circ P$, where $P=P_0$,
$E_0$ is the identity map, $c=min\,\{-a,b\}$, \mbox{$\delta\in
(0,c/8)$,} and \mbox{$E_{\lambda}:[a,b] \to [a,b]$} is a $C^r$
diffeomorphism satisfying the following conditions:
\begin{eqnarray}
E_\lambda(x)-x &=& \lambda\delta, \quad x\in [-4\delta,
 4\delta]\label{eq1},\\
E_\lambda(x)-x&\le&\lambda\delta,\quad x\in [a,b].\label{eq2}
  \end{eqnarray}
\end{lem}
\proof By Theorem \ref{fbt}, there exist $\epsilon>0$ arbitrarily small, a flow box
\mbox{$B=B([a,b],\epsilon)\subset V_0$}, and a
$C^{r+1}-$rectifying diffeomorphism $h:B\to [-\epsilon,0]\times
[a,b]$. Let $\phi_1:[-\epsilon,0]\to
[0,1]$ and $\phi_2:[a,b]\to [0,1]$ be smooth functions such that \;
$(\phi_1)^{-1}(1) = [-8\epsilon/10,-2\epsilon/10]$,\;
$(\phi_1)^{-1}(0) = [-\epsilon,0]\setminus [-9\epsilon/10,
-\epsilon/10]$, \;
%$\phi_1(x-\epsilon/2)\equiv\phi_1(-x-\epsilon/2)$,
$(\phi_2)^{-1}(1) = [-6\delta,6\delta]$, \; $(\phi_2)^{-1}(0) =
[a,b]\setminus [-7\delta,7\delta]$.
%, and $\phi_2(x)\equiv\phi_2(-x)$.
Let $Y_0:[-\epsilon,0]\times [a,b]\to \R^2$ be the smooth vector
field which at each $(x,y)\in [-\epsilon,0]\times [a,b]$ takes the
value:
$$Y_0(x,y)=(1,0)+\eta\phi_1(x)\phi_2(y)(0,\delta),$$
where $\eta>0$ is a positive constant such that the positive
semitrajectory $\gamma_{(-\epsilon,-4\delta)}^+$ of $Y_0$ starting
at $(-\epsilon,-4\delta)$ intersects $\{0\}\times [a,b]$ at the
point $(0,-3\delta)$. By construction, for each $y\in
[-4\delta,4\delta]$, the positive semitrajectory
$\gamma_{(-\epsilon,y)}^+$ of $Y_0$ starting at $(-\epsilon,y)$ is
an upward shift of ${\gamma}_{(-\epsilon,-4\delta)}^+$ and so
intersects $\{0\}\times [a,b]$ at $(0,y+\delta)$. Define
$Y\in\mathfrak{X}_H^r(M)$ to be a vector field supported in $B$ such
that $Y\vert_{B}=(h^{-1})_*(Y_0)$. Accordingly,
$$(X+\lambda Y)\vert_B=(h^{-1})_*((1,0)+\lambda Y_0).$$
Recall that by Theorem~\ref{fbt}, the map  $h$ takes isometrically
$[a,b]$ onto $\{0\}\times [a,b].$ By construction, the
one--parameter family of vector fields $X+\lambda Y$ has all the
required properties.
\endproof

%\begin{lem}\label{lem9}
%Let $[a,b]$ be a closed subinterval which is a transversal section
%to \linebreak $X\in\mathfrak{X}^\infty(M)$. Let $\varepsilon > 0$ be
%such that $B=\{X_t (p) : p \in [a,b], \; t \in [-\varepsilon,0]\}$
%is a flow box. Let $Y\in\mathcal{A}(B,X),$ $\lambda\in[0,1]$ and
%$T_\lambda :[a,b]\to [a,b]$ be the forward Poincar\'e Map induced by
%$X+\lambda Y$ on $[a,b]$. Then
%\item [$(a)$] If $X$ is smooth, $a<0<b$ and $c=min\{-a,b\},$ the vector field
%$Y$ can be taken so that for every $\lambda \in [0,1]$, there
%  exists a $C^r$ diffeomorphism $E_{\lambda}:[a,b] \to [a,b]$
%  such that $T_\lambda=E_{\lambda} \circ T$, where $E_0$ is the
%  identity map, $$0<\sup \{|E_1(x)-x| : x \in [a,b]\}= \delta\in
%  (0,c/8),$$ and  for all $\, x\in [-4\delta,
% 4\delta],\, $ $\,
%E_\lambda(x)-x = \lambda\delta;$
%\end{lem}

\begin{teo}\label{tt} Let $X\in\mathfrak{X}_H^r(M)$, $\sigma$  be
a non--trivial recurrent unstable saddle--separatrix,
\mbox{$\Sigma=[a,b]$} be a $\sigma$-adapted transversal
segment to $X$, $B=B([a,b],\epsilon)$ be a flow box of $X$ and
 \mbox{$Y\in
\mathcal{A}(B,X)$.} If $q\in [a,b]$ is the first intersection of
$\sigma$ with $[a,b]$ then either of the following alternatives
happens:

  \mm $(a)$ for some $\lambda\in[0,1],$  $\; [a,b]$ intersects a
  saddle--connection of $\,X+\lambda Y$ or,

  \mm $(b)$ for every $(\lambda,n)\in [0,1]\times \N$, the point $q$ belongs to
   ${\rm dom}\, (P_\lambda^n)$
   %$($resp. $q$ belongs to
  %$\mbox{Dom}\,((T_\lambda)^{-n}))$
   and $P_\lambda^n(q)$
   %(resp.
  %$(T_\lambda)^{-n}(q)$)
   depends continuously on $\lambda$. In
  this case, for each $\lambda\in [0,1]$, the sequence $\{P^n_{\lambda}(q)\}_{n\in\N}$
  %(resp.
  %$\{T^{-n}_{\lambda} (p)\}$)
   accumulates in a point of $[a,b]$
  belonging, with respect to $X+\lambda Y,$ to either a closed
  trajectory or to
   a non-trivial recurrent trajectory, where $P_\lambda:[a,b]\to
  [a,b]$ denotes the forward Poincar\'e map induced by $X+\lambda Y$.
\end{teo}
\proof Assume that $(a)$ does not happen. Let us prove that then
$(b)$ occurs. Firstly we have to show that for every $(\lambda,n)\in
[0,1]\times\N$, the point $q$ belongs to ${\rm dom}\,
(P_\lambda^n)$. Suppose that this does not happen. So for some
\mbox{$(\lambda_1,n_1)\in (0,1]\times\N-\{0\}$,} we have that $q\in
{\rm dom}\,(P_\lambda^{n_1-1})$ for all \mbox{${\lambda}\in [0,1]$,}
and $q\not\in {\rm dom}\,(P_{\lambda_1}^{n_1})$. Hence, we have that
$P_{\lambda_1}^{n_1-1}(q)$ does not belong to ${\rm
dom}\,(P_{\lambda_1})={\rm dom}\,(P_{0})$ whereas $P_0^{n_1-1}(q)\in
{\rm dom}\,(P_0)$. By construction, $P_\lambda^{n_1-1}(q)$ depends
continuously \mbox{on $\lambda$,} and so for some $\lambda_2\in
[0,\lambda_1]$, $P_{\lambda_2}^{n_1-1}(q)$ intersects the boundary
of ${\rm dom}\,(P_0)$. By Lemma \ref{lem8}, $X+\lambda_2 Y$ has a
saddle--connection intersecting $[a,b]$, which contradicts the
initial assumption. Therefore, the first part of $(b)$ is proved.
The second part of $(b)$ follows from Lemma \ref{td2} since the
existence of an attracting graph intersecting $[a,b]$ would \mbox{imply
(a).}
\endproof

In the proof of next lemma we shall use the fact that a transversal
segment $\Sigma=[a,b]$ to $X\in\mathfrak{X}_H^r(M)$
 may also be represented by $[a+s,b+s]$, for any $s\in\R$. Henceforth, if $A$ is a subset of $M$ then $\overline{A}$ will
denote its topological closure.

\begin{lem}\label{lem10} Let $X\in\mathfrak{X}_H^r(M)$, $r\ge 2$, be
smooth in a neighborhood $V_0$ of a non--trivial \mbox{recurrent}
point $p\in M$. Assume that $X$ has the infinitesimal contraction
property at $p$. Given a neighborhood $V$ of $p$, there exist a flow
box $B\subset V$ and \mbox{$Y\in \mathcal{A}(B,X)$}
\mbox{arbitrarily} $C^r-$close to the \mbox{zero--vector--field}
such that for some \mbox{$\lambda\in [0,1]$,} \mbox{$X+\lambda Y$}
has a saddle--connection \mbox{meeting $B$.}
\end{lem}
\proof By Lemma \ref{Peixoto2}, there exist non--trivial recurrent
saddle--separatrices $\sigma_1$, $\sigma_2$ such that
\mbox{$\omega(\sigma_2)\cap\alpha(\sigma_1)=\overline{\gamma_p}$.}\,
Let {$\Sigma_1=[a_1,b_1]\subset V_0\cap V$} be a transversal
 segment to $X$
passing through $p$ such that $P_{\Sigma_1}$ is an infinitesimal
$\kappa-$contraction for some $0<\kappa<0.1$. By \mbox{Lemma
\ref{lem7},} there exists a $\sigma_2$--adapted subsegment
$\Sigma=[a,b]\subset [a_1,b_1]$. Let $\overline{p}\in (a,b)$ be the
first intersection of $\sigma_1$ with $(a,b)$. Accordingly,
$\overline{p}$ is a non--trivial recurrent point. Modulo shifting
the interval $[a_1,b_1]$, we may assume that $a<0<b$ and
$\overline{p}=0$. Let $B=B([a,b],\epsilon)\subset V_0\cap V$ be a
flow box for some $\epsilon>0$. By Lemma \ref{lem9}, there exists
$Y\in\mathcal{A}(B,X)$ arbitrarily $C^r-$close to the
zero--vector--field such that the forward Poincar\'e Map
$P_\lambda=E_\lambda\circ P$ induced by $X+\lambda Y$ on $[a,b]$ has
the properties (\ref{eq1}) and (\ref{eq2}). We shall consider only
the case in which $0$ is an accumulation point of $\sigma_2\cap
[a,0)$. Let $q\in\sigma_2\cap [a,b]$ be the first intersection of
$\sigma_2$ with $[a,b]$.

Suppose by contradiction that, for all $\lambda \in [0,1],$ $\;
X+\lambda Y$ has no saddle--connections. Then by Theorem \ref{tt},
for all $(\lambda,n)\in [0,1]\times \N$, the point $q$ belongs to
   ${\rm dom}\, (P_\lambda^n)$
   %$($resp. $q$ belongs to
  %$\mbox{Dom}\,((T_\lambda)^{-n}))$
   and $P_\lambda^n(q)$
   %(resp.
  %$(T_\lambda)^{-n}(q)$)
   depends continuously on $\lambda$. By $(\ref{eq2})$ of Lemma \ref{lem9} and by proceeding inductively, we
may see that, for all integer $n \geq 1$,
$$
|P \circ (E_{\lambda} \circ P)^{n-1} (q)-P^n(q)|\; \leq  \; \kappa
\delta (1+\kappa+\cdots+\kappa^{n-2}) \leq  \;
\frac{\kappa\delta}{1-\kappa}\,.
$$
%and
%$$
%|(E_{\lambda} \circ P)^{n} (q)-P^n(q)|\; \leq  \; \delta + \kappa
%\delta (1+\kappa+\cdots+\kappa^{n-2}) \leq  \; \delta +
%\frac{\kappa\delta}{1-\kappa}\,.
%$$

As $0$ is an accumulation point of $\sigma_2 \cap [a,0)$ there
exists $N \in \N$ such that \linebreak\mbox{$P^N(q) \in
[-\kappa\delta,0]$.} Therefore,
$$
P \circ (E_1 \circ P)^{N-1}(q) \, \geq \, P^N (q) - \frac{\kappa
\delta}{1-\kappa} \, \geq \, -\kappa \delta - \frac{\kappa
\delta}{1-\kappa} \, \geq \, -3\kappa \delta.
$$

Hence, by $(\ref{eq1})$ of  Lemma \ref{lem9} and by the fact that
$0<\kappa<0.1$,
$$
(E_1 \circ P)^{N}(q) \, = \, E_1 \circ (P \circ (E_1 \circ P)^{N-1})
(q) \, = \, P \circ (E_1 \circ P)^{N-1} (q)+\delta \, \geq \,
-3\kappa \delta + \delta > 0\,.
$$

This implies that there exists $\lambda \in [0,1]$ such that
$P_\lambda^{N}(q)=(E_{\lambda} \circ P)^N (q)=0$ (see (b) of Theorem
\ref{tt}). That is, $X+\lambda Y$ has a saddle--connection passing
through $0$. This contradiction proves the lemma. \endproof

\bb

\begin{teo}\label{saddle-connection} Suppose that $\mathfrak{X}_H^r(M)$, $r\ge 2$, has
the infinitesimal contraction property at a non--trivial recurrent
point $p$. Then, given neighborhoods $V$ of $p$ in $M$ and
$\mathcal{V}$ of $X$ in $\mathfrak{X}^r(M)$, there exist $Z\in \mathcal{V}$, with
$Z-X$ supported in $V$, having either a periodic trajectory
meeting $V$ or a saddle--connection meeting $V$.
\end{teo}

\proof Let be given neighborhoods $V$ of $p$ in $M$ and
$\mathcal{V}$ of $X$ in $\mathfrak{X}_H^r(M)$. By \mbox{Lemma
\ref{pro5},} there exist a a flow box $B_0\subset V$ and a
neighborhood $\mathcal{V}_0\subset\mathcal{V}$ of $X$ in
$\mathfrak{X}_H^r(M)$ such that every $Z\in\mathcal{V}_0$, with
$Z-X$ supported in $B_0$, has the infinitesimal contraction property
at $B_0$. By the proof of Lemma \ref{pro5} and by Lemma \ref{lem7},
we may assume that $B_0=B(\Sigma,\epsilon)$, where $\Sigma$ is a
$\sigma-$adapted transversal segment to $X$ for some non--trivial
recurrent unstable saddle--separatrix $\sigma$. By shrinking
$\mathcal{V}_0$ if necessary, we may assume that
 for every $Z\in\mathcal{V}_0$ with
$Z-X$ supported in $B_0$ we have that $Z-X\in\mathcal{A}(B,X)$.
Suppose, by contradiction, that  every vector field in
$\mathcal{V}_0$ with $Z-X$ supported in $B_0$ has neither periodic
trajectories meeting $B_0$ nor saddle--connections meeting $B_0$. We
claim that, under these assumptions,  every $Z\in \mathcal{V}_0$
with $Z-X$ supported in $B_0$ has a non-trivial recurrent point in
the interior of $B_0$. Indeed, by taking $\lambda=1$ in $(b)$ of
Theorem \ref{tt}, we get that every $Z=X+(Z-X)\in\mathcal{V}_0$ with
$Z-X$ supported in $B_0$ has a non--trivial recurrent point
intersecting the boundary of $B_0$. Since $B_0$ is still a flow box
of $Z$, we have that the interior of $B_0$ has non--trivial
recurrent points of $Z$. This proves the claim. Now let
$Z_1\in\mathcal{V}_0$ be a $C^r$ vector field which is smooth in
$B_0$ and is such that $Z_1-X$ supported in $B_0$. By the claim,
$Z_1$ has a non--trivial recurrent point $p_1$ in the interior of
$B_0$, and $Z_1$ has the infinitesimal contraction property at
$B_0$. By \mbox{Lemma \ref{lem10}}, there exist a flow box $B\subset
V$ and $Z_2\in\mathcal{V}_0$, with $Z_2-X$ supported in $B$, having
a saddle--connection meeting $B$. This contradiction finishes the
proof.
\endproof

\section{$C^r-$Closing Results}

An {\it interval exchange transformation} or simply an {\it iet} is
an injective map $E:\R/\Z\to\R/\Z$ of the unit circle,
differentiable everywhere except possibly at finitely many points
and such that for all $x\in {\rm dom}\,(E)$ (its domain of
definition), $\vert DE(x)\vert=1$. The trajectory of $E$ passing
through $x\in\R/\Z$ is the set
$O(x)=\{E^n(x):n\in\Z\,\,\,\text{and}\,\,\, x\in{\rm dom}\,(E^n)\}$.
We say that $E$ is {\it minimal} if $O(x)$ is dense in $\R/\Z$ for
every $x\in\R/\Z$. Given a transversal circle $C$ to
$X\in\mathfrak{X}_H^r(M)$, we say that the forward Poincar\'e Map
$P:C\to C$ is {\it topologically semiconjugate} to an iet
$E:\R/\Z\to\R/\Z$ if there is a monotone continuous map
$h:C\to\R/\Z$ of degree one such that $E\circ h(x)=h\circ P(x)$ for
all $x\in{\rm dom}\,(P)$.

We shall need the following structure theorem due to Gutierrez
\cite{G2}. We should remark that in this theorem below, the item
$(d)$ although not explicitly stated in \cite{G2} follows from the
proof given therein and from the fact that $X$ has finitely many
singularities.

\begin{teo}\label{gst} Let
$X\in\mathfrak{X}_H^r(M)$. The topological closure of the
non--trivial recurrent trajectories of $X$ determines finitely many
quasiminimal sets $N_1$,$N_2$\ldots, $N_m$. For each $1\le i\le m$,
there exists a transversal circle $C_i$ to $X$ intersecting every
non--trivial recurrent trajectory of $X\vert_{N_i}$ such that if
$P_i:C_i\to C_i$ is the forward Poincar\'e Map induced by $X$ on
$C_i$ then:
\begin{itemize}
\item [$(a)$] Either $N_i\cap C_i=C_i$ or
$N_i\cap C_i$ is a Cantor set;\vspace{0.2cm} \item[$(b)$] $N_j\cap
C_i=\emptyset$, for all $j\in\{1,2\ldots,i-1,i+1,\ldots,m\}$;
\item [$(c)$] $P_i$ is topologically semiconjugate to a minimal
interval exchange transformation \linebreak $E_i:\R/\Z\to\R/\Z$; %in the
%case in which $\Lambda={\rm cl}\,(\gamma)\cap C$ is a Cantor set;
%Furthermore, the complement of the Cantor set $\Lambda$ consists of
%countable many disjoint open intervals, called adjacent intervals,
%each of which is carried by the semiconjugacy $h$ into a point; and
%the restriction of $h$ to the complement of the closures of all the
%adjacent intervals is a strictly monotone map;
\item [$(d)$] For each $q\in C_i$,
$\gamma_q\cap C_i$ is an infinite set.
\end{itemize}
We call the circle $C_i$ a special transverse circle for $N_i$.
\end{teo}

\begin{cor}\label{tend} Let $X\in\mathfrak{X}_H^r(M)$ and let $N$ be
a quasiminimal set. Given a transversal segment $\Sigma_1$ passing
through a non--trivial recurrent point $p\in N$, there exists a
subsegment $\Sigma$ of $\Sigma_1$ passing through $p$ such that if
$z\in\Sigma$ then either $\alpha(z)=N$ or $\omega(z)=N$. In
particular, either $z\in\cap_{n=1}^\infty\,{\rm dom}\,(P^n)$ or
$z\in\cap_{n=1}^\infty\,{\rm dom}\,(P^{-n})$, where
$P:\Sigma\to\Sigma$ is the forward Poincar\'e Map induced by $X$.
\end{cor}
\proof Let $C$ be a special transversal circle for $N$. There exist
a subsegment $\Sigma$ of $\Sigma_1$ passing through $p$ and a
subsegment $\Gamma$ of $C$ such that the forward Poincar\'e Map
$T:\Sigma\to\Gamma$ induced by $X$ is a diffeomorphism. Since $C$ is
free of finite trajectories (by $(d)$ of Theorem \ref{gst}), so is
$\Sigma$. Hence, by Lemma \ref{td2}, either $\alpha(z)$ or
$\omega(z)$ is a quasiminimal set, which by $(b)$ of Theorem
\ref{gst}, has to be $N$.
\endproof

\begin{pro}\label{clcp} Suppose that $X\in\mathfrak{X}_H^r(M)$ has the infinitesimal
contraction property at a non--trivial recurrent point $p\in M$.
There exists an arbitrarily small flow box $B_0$ ending at $p$ and an arbitrarily
small neighborhood $\mathcal{V}_0$ of $X$ in $\mathfrak{X}_H^r(M)$
such that either:
\begin{itemize}
\item [$(i)$] some $Z\in\mathcal{V}_0$ with $Z-X$ supported in $B_0$ has a periodic trajectory
meeting $B_0$ or,
\item[$(ii)$] every $Z\in\mathcal{V}_0$ with $Z-X$ supported in $B_0$ has a non--trivial recurrent point in the interior
of $B_0$.
\end{itemize}
\end{pro}
\proof By Corollary \ref{tend}, given a transversal segment $\Sigma_1$ to $X$
passing through $p$, there exists a subsegment $\Sigma$ of $\Sigma_1$ passing through $p$ such that for every
$z\in\Sigma$, either $\alpha(z)=N$ or $\omega(z)=N$, where
$N=\overline{\gamma_p}$. By taking a subsegment of $\Sigma$ if
necessary, we may assume that the forward Poincar\'e Map
$P:\Sigma\to\Sigma$ induced by $X$ is an infinitesimal
$\kappa-$contraction for some $\kappa\in (0,1)$.

We claim that
$z\in\Sigma\setminus \cap_{n=1}^\infty \,{\rm dom}\,(P^n)$ if and
only if $\omega(z)$ is a saddle--point. Indeed, if
$z\in\Sigma\setminus \cap_{n=1}^\infty \,{\rm dom}\,(P^n)$ then
there exists $m\in\N$ such that $z\in\,{\rm dom}\,(P^m)$ but
$z\not\in {\rm dom}\,(P^{m+1})$. Hence, $P^m(z)\not\in {\rm
dom}\,(P)$ and by \mbox{Lemma \ref{lem8},} either $\omega(z)$ is a
saddle--point or $P^m(z)$ belongs to the open set $\Sigma\setminus
{\rm dom}\,(P)$. In this last case, there exists a subsegment
$I\subset\Sigma$ containing $z$ such that $P^m(I)\subset \Sigma\setminus {\rm
dom}\,(P)$ and $I\subset \cap_{n=1}^\infty \,{\rm dom}\,(P^{-n})$
(by the first part of this proof). Of course, this is impossible
since $P^{-1}$ has a uniformly expanding behaviour and $\Sigma$ has
finite length. This proves the claim.

In particular, we have that ${\rm dom}\,(P)$ is the whole
transversal segment $\Sigma$ but finitely many points. Let
$B_0=B(\Sigma,\epsilon)$ be a flow box and let
$\mathcal{V}_0\subset\mathfrak{X}_H^r(M)$ be a neighborhood of $X$
such that if $Z\in \mathcal{V}_0$ and $Z-X$ is supported in $B_0$
then $B_0$ is still a flow box of $Z$ and so \mbox{${\rm
dom}\,(P_Z)= {\rm dom}\,(P)$}, where $P_Z$ is the forward Poincar\'e
Map induced by $Z$ on $\Sigma$. Hence, for every
${Z\in\mathcal{V}_0}$ such that $Z-X$ is supported in $B_0$, ${\rm
dom}\,(P_Z)$ is the whole transversal segment but finitely many
points whose positive trajectories go directly to saddle--points.
Since there are only finitely many saddle--points, we have
that for each $Z\in\mathcal{V}_0$ such that $Z-X$ is supported in
$B_0$, there exists a countable subset $D$ of $\Sigma$ such that for
every $z\in\Sigma\setminus D$ the positive semitrajectory of $Z$
starting at $z$ intersects $\Sigma$ infinitely many times. By \mbox{Lemma
\ref{td2},} $\omega(z)$ is either a recurrent trajectory intersecting
$B_0$ or an attracting graph intersecting $B_0$. In the second case,
an arbitrarily small $C^r-$perturbation of $Z$ supported in $B_0$
yields a vector field $\widetilde{Z}\in\mathcal{V}_0$
 having a
periodic trajectory meeting $B_0$. \endproof

\begin{teo}[Localized $C^r-$Closing Lemma]\label{mt1} Suppose that $X\in\mathfrak{X}_H^r(M)$, $r\ge 2$, has the
contraction property at a non--trivial recurrent point $p\in M$ of
$X$. Given neighborhoods $V$ of $p$ in $M$ and $\mathcal{V}$ of $X$
in $\mathfrak{X}_H^r(M)$, there exists $Y\in\mathcal{V}$, with $Y-X$
supported in $V$, such that $Y$ has a periodic trajectory meeting
$V$.
\end{teo}

\proof Assume by contradiction that no vector field
$Y\in\mathcal{V}$ with $Z-X$ supported in $V$ has a periodic
trajectory meeting $V$. By Proposition \ref{clcp} and by Lemma
\ref{pro5}, there exist a flow box $B_0\subset V$ and a neibhborhood
$\mathcal{V}_0\subset\mathcal{V}$ of $X$ such that every
$Z\in\mathcal{V}_0$ with $Z-X$ supported in $B_0$ has the
infinitesimal contraction property at $B_0$ and a non--trivial
recurrent point in ${\rm int}\,(B_0)$, the interior of $B_0$. Note
that every vector field  $Z\in\mathcal{V}_0$ with $Z-X$ supported in
$B_0$ has at most $4N_s$ saddle--connections, where $N_s$ is the
number of saddle--points of $X$. Therefore, the proof will be
finished if we construct a sequence $\{Z_n\}_{n=0}^{4N_s+1}$ of
vector fields in $\mathcal{V}_0$ such that for each $n\in\N$,
$Z_n-X$ is supported in $B_0$ and $Z_{n+1}$ has one more
saddle--connection than $Z_n$. Let us proceed with such a
construction. Let $p_0\in {\rm int}\,(B_0)$ be a non--trivial
recurrent point of $Z_0=X$. By Theorem \ref{saddle-connection},
there exist an open set $V_1\subset B_0$ and $Z_1\in\mathcal{V}_0$
with $Z_1-X$ supported in $V_1$ having a saddle--connection
$\sigma_1$ meeting $V_1$. By the above, $Z_1$ has also a
non--trivial recurrent point $p_1\in{\rm int}\,(B_0)$. Now we may
repeat the reasoning. By \mbox{Theorem \ref{saddle-connection},}
there exist an open set $V_2\subset B_0\setminus\sigma_1$ and
$Z_2\in\mathcal{V}_0$ with $Z_2-X$ supported in $V_2$ having a
saddle--connection $\sigma_2$ meeting $V_2$ (and a
saddle--connection $\sigma_1$ meeting $V_1$). Moreover, $Z_2$ has a
non--trivial recurrent point $p_2\in{\rm int}\,(B_0)$. By proceeding
by induction, we shall obtain a vector field
$Z_{4N_s+1}\in\mathcal{V}_0$ with $Z_{4N_s+1}-X$ supported in $B_0$
having at least $4N_s+1$ saddle--connections, which is a
contradiction.
\endproof

\noindent{\rm\bf Theorem A.} {\it Suppose that
$X\in\mathfrak{X}_H^r(M)$, $r\ge 2$, has the contraction property at
a {quasimi\-nimal} set $N$. For each $p\in N$, there exists
\mbox{$Y\in\mathfrak{X}_H^r(M)$} arbitrarily $C^r-$close to $X$
having a periodic trajectory passing through $p$.} \proof That
localized $C^r-$closing (Theorem \ref{mt1}) implies $C^r-$closing
(Theorem A) is an elementary fact.
\endproof

\section{Transverse Measures}

Let $N$ be a quasiminimal set of $X\in\mathfrak{X}_H^r(M)$, $\Sigma$
be a transversal segment to $X$ such that
$\Sigma\setminus\partial\Sigma$ intersects $N$ and
$P:\Sigma\to\Sigma$ be the forward Poincar\'e Map induced by $X$. We
may consider $\Sigma$ as a Borel measurable space
$(\Sigma,\mathscr{B})$, where $\mathscr{B}$ is the Borel
$\sigma-$algebra on $\Sigma$. We say that a Borel probability
measure is {\it non--atomic} if it assigns measure zero to every
one--point--set. A {\it transverse measure on $\Sigma$} is a
non--atomic $P-$invariant Borel probability measure which is
supported in $N\cap \Sigma$. A transverse measure $\nu$ is called
{\it ergodic} if whenever $P^{-1}(B)=B$ for some Borel set
$B\in\mathscr{B}$ then either $\nu(B)=0$ or $\nu(B)=1$. We let
$\mathscr{M}(\Sigma)$ denote the set of Borel probability measures
on $\Sigma$ and we let $\mathscr{M}_P(\Sigma)$ denote the subset of
$\mathscr{M}(\Sigma)$ formed by the $P-$invariant Borel probability
measures. We say that $P$ is {\it uniquely ergodic} if
$\mathscr{M}_P(\Sigma)$ is a unitary set. A set $W\subset\Sigma$ is
called a {\it a total measure set} if $\nu(W)=1$ for every
$\nu\in\mathscr{M}_P(\Sigma)$. Concerning the existence of
transverse measures, we have the following result.

\begin{teo}\label{invarm} Let $N$ be a quasiminimal set of
$X\in\mathfrak{X}_H^r(M)$ and let $\Sigma_1$ be a compact transversal segment
to $X$ passing through a non--trivial recurrent point $p\in N$. There \mbox{exist} a subsegment $\Sigma\subset\Sigma_1$ passing through $p$ and finitely many ergodic transverse measures $\nu_1,\ldots,\nu_s\in {\mathscr{M}_P(\Sigma)}$
such that every $\nu\in\mathscr{M}_P(\Sigma)$ can be written in the form $\nu=\sum_{i=1}^s\lambda_i\nu_i$,
where $\lambda_i\ge 0$ for all $1\le i\le s$, and $\sum_{i=1}^s\lambda_i=1$.
\end{teo}
\proof The proof may be split into two parts. The first part of the
proof -- that every small subsegment of $\Sigma_1$ passing through
$p$ can be endowed with a transverse measure -- can be found in
Katok \cite{Kat} and \mbox{Gutierrez \cite{G9}.} To prove the second
part, let $C$ be a special transversal circle to $X$ passing through
$\gamma_p$ as in the \mbox{Theorem \ref{gst}.} There exist
subsegments $\Sigma\subset\Sigma_1$ containing $p$ and
$\Gamma\subset C$ such that the forward Poincar\'e Map
$T:{\Sigma}\to\Gamma$ induced by $X$ is a diffeomorphism. We claim
that $\mathscr{M}_P({\Sigma})$ is made up of transverse measures,
where $P:\Sigma\to\Sigma$ is the forward Poincar\'e Map induced by
$X$. Indeed, by $(d)$ of Theorem \ref{gst}, \mbox{$\Sigma$ is} free
of periodic points. By Poincar\'e \mbox{Recurrence} Theorem, the set
of non--trivial recurrent points in $\Sigma$ is a total measure set.
\mbox{By $(b)$} of Theorem \ref{gst}, all these non--trivial
recurrent points belong to the same quasiminimal set. This proves
the claim. Now, every (ergodic) transverse measure on $\Sigma$
 corresponds, via the diffeomorphism $T$, to a (ergodic)
transverse measure on $C$. By $(c)$ of \mbox{Theorem \ref{gst},}
every (ergodic) transverse measure on $C$ corresponds to a (ergodic)
Borel probability measure on $\R/\Z$ invariant by a minimal interval
exchange transformation \mbox{$E:\R/\Z\to\R/\Z$.} By a result of
\mbox{Keane \cite{Kea},} which also holds for interval exchange
transformations with flips \cite{CoFoSi}, there exist only finitely
many ergodic Borel probability measures invariant by $E$. Each of
such $E-$invariant Borel probability measures on $\R/\Z$ is
associated to exactly one ergodic transverse measure in
$\mathscr{M}_P({\Sigma})$. Now the rest of the proof follows from
the fact that $\mathscr{M}_P(\Sigma)$ is the convex hull of its
ergodic measures.
\endproof

Let $P:\Sigma\to\Sigma$ be the forward Poincar\'e Map induced by $X$
on a transversal segment $\Sigma$ to $X\in\mathfrak{X}_H^r(M)$. By
\mbox{Lemma \ref{lem8},} the domain of $P$ is the union of finitely
many open, pairwise disjoint subintervals of $\Sigma$: ${\rm
dom}\,(P)=\cup_{i=1}^m J_i$. We say that the {\it lateral limits of
$\vert DP\vert$
 exist} if for every $1\le i\le m$ and for every $p\in \partial J_i$, the lateral limit
$\ell=\lim_{x\to p\atop x\in J_i} \vert DP(x)\vert$ exists as a
point of $[0,+\infty]$.

Henceforth, till the end of this paper, we shall assume that $N$ is
a quasiminimal set, $\Sigma$ is a transversal segment to $X$ such
that $\Sigma\setminus\partial\Sigma$ intersects $N$ and
$P:\Sigma\to\Sigma$ is the forward Poincar\'e Map induced by $X$ on
$\Sigma$. We shall assume that $\Sigma$ is so small that the forward
Poincar\'e Map $T:\Sigma\to T(\Sigma)\subset C$ induced by $X$ is a
diffeomorphism, where $C$ is a special transversal circle for $N$,
and that
 $P$ has the following properties:
\begin{itemize}
 \item [$(P1)$] $\vert DP\vert$ is bounded from above;
 \item [$(P2)$] The lateral limits of $\vert DP\vert$ exist.
\end{itemize}

\medskip

\begin{defi}[Almost--integrable function] We say that $\log\vert DP\vert$
is $\nu$--{\it almost--integrable} if
\[{\rm min}\,\Big\{\int {{\log}^+\vert DP\vert} {\,\mathrm d}\nu,\int {{\log}^-\vert DP\vert}{\,\mathrm d}\nu
\Big\}<\infty,
\]
where
\[ {\log^+\vert DP(x)\vert}={\rm max}\,\{\log\vert
DP(x)\vert,0\}, \quad {\log^-\vert DP(x)\vert}={\rm
max}\,\{-\log\vert DP(x)\vert,0\},\] and
$\nu\in\mathscr{M}(\Sigma)$. In this case we define
\[\int \log\vert DP\vert {\,\mathrm d}\nu=\int {\log^+\vert DP\vert}{\,\mathrm d}\nu-\int {\log^-\vert DP\vert}{\,\mathrm d}\nu,
\]
which is a well defined value of the subinterval $[-\infty,\infty)$
of the extended real line $[-\infty,\infty]$.
\end{defi}

\begin{lem}\label{cext} Suppose that there exists $K\in\R$ such that $\int\log\vert DP\vert\,{\rm d}\nu<K$
 for all\, \mbox{$\nu\in\mathscr{M}_P(\Sigma)$.} Then there exists a continuous function $\phi:\Sigma\to\R$
everywhere defined, with \mbox{$\log\vert DP(x)\vert<\phi(x)$} for
all $x\in {\rm dom}\,(P)\setminus P^{-1}(\partial\Sigma)$, such that
$\int\phi\,{\rm d}\nu<K$ for all\, $\nu\in\mathscr{M}_P(\Sigma)$.
\end{lem}
\proof By reasoning as in Theorem \ref{invarm}, since $\Sigma$ is
disjoint of periodic trajectories, we may show that
$\mathscr{M}_P(\Sigma)$ is the convex hull of finitely many ergodic
(non--atomic) transverse measures $\nu_1$, $\ldots$, $\nu_s$. It
follows from $(P1)$ and $(P2)$ that there exists a continuous
function $\overline{\phi}:\overline{{\rm dom}\,(P)}\to\R$ such that
$\int\overline{\phi}\,{\rm d}\nu_i<K$, for all $1\le i\le s$, and
$\log\vert DP(x)\vert<\overline{\phi}(x)$ for all $x\in{\rm
dom}\,(P) \setminus P^{-1}(\partial\Sigma)$. Hence,
$\int\overline{\phi}\,{\rm d}\nu<K$ for all
$\nu\in\mathscr{M}_P(\Sigma)$. Now we may take $\phi$ to be any
continuous extension of $\overline{\phi}$ to $\Sigma$. Since every
$\nu\in\mathscr{M}_P(\Sigma)$ is supported in
$N\cap\Sigma\subset\overline{{\rm dom}\,(P)}$, we have that
$\int\phi\,{\rm d}\nu=\int\overline{\phi}\,{\rm d}\nu<K$ for all
$\nu\in\mathscr{M}_P(\Sigma)$.
\endproof

\begin{lem}\label{abc} The following statements are equivalent:
\begin{itemize}
 \item [$(a)$] $\liminf_{n\to\infty}\frac1n\log\vert DP(x)\vert<0$ for all $x$ in a total measure set;
 \item [$(b)$] $\int\log\vert DP\vert\,{\rm d}\nu<-c$ for some $c>0$ and for all $\nu\in\mathscr{M}_P(\Sigma)$;
 \item [$(c)$] $\liminf_{n\to\infty}\frac1n\log\vert DP(x)\vert<-c$ for some $c>0$ and for all $x$ in a total measure set;
\end{itemize}
\end{lem}
\proof Let us show that $(a)$ implies $(b)$. By $(P1)$, $\log\,\vert
DP\vert$ is $\nu-$almost integrable with respect to each
$\nu\in\mathscr{M}(\Sigma)$. Hence, there exists $K\in\R$ such that
 \mbox{$\int\log\vert DP\vert\,{\rm d}\nu_i<K$,} for all \mbox{$1\le i\le
 s$,}
where $\{\nu_i\}_{i=1}^s$ are the ergodic transverse measures in
$\mathscr{M}_P(\Sigma)$. So either \mbox{$\int\log\vert
DP\vert\,{\rm d}\nu_i=-\infty$} for all $i=1,\ldots,s$ (and we are
done) or there exists a non--empty subset $\Lambda$ of
$\{1,2,\ldots,m\}$ such that $\log\vert DP\vert$ is
$\nu_i-$integrable for all $i\in\Lambda$. In this case, $(a)$ and
Birkhoff Ergodic Theorem yields that
$$\int\log\vert DP\vert\,{\rm
d}\nu_i=\lim_{n\to\infty}\frac1n\log\vert
DP(x)\vert=\liminf_{n\to\infty}\frac1n\log\vert DP(x)\vert=-c_i<0$$
for some $x$ in a $\nu_i$--full measure set. Now take
$c=\min\,\{c_i:i\in\Lambda\}$. A similar reasoning shows that $(b)$
implies $(c)$. This finishes the proof.
\endproof

\begin{lem}\label{Par}  Let $\{\mu_j\}_{j\in\N}$ be a sequence of Borel probability measures in $\mathscr{M}(\Sigma)$
weakly$^*$ converging to $\mu\in \mathscr{M}(\Sigma)$. The following hold:
\begin{itemize}
 \item [$(a)$] $\mu(B)=\lim_{j\to\infty}\mu_j({B})$ for every Borel set $B\in\mathscr{B}$ such that
$\mu(\partial B)=0$, where $\partial B$ denote the topological boundary of $B$;
\item [$(b)$] $\mu(J)=\lim_{j\to\infty}\mu_j({J})$ for every open subinterval $J$ of $\Sigma$ such that
$\mu(\partial J\setminus\partial\Sigma)=0$.
\end{itemize}
\proof The item $(a)$ is a standard theorem from measure theory (see
\cite[Theorem 6.1]{Par}). Let us prove $(b)$. Let $J$ be an open
subinteval of $\Sigma$. If $\partial J\cap\partial\Sigma=\emptyset$
then $\mu(\partial J)=\mu(\partial J\setminus\partial\Sigma)=0$ and
the result follows from $(a)$. If $J=\Sigma$ then the indicator
function $\chi_J$ is continuous and so the result follows
immediately from the weak$^*$ convergence of $\{\mu_j\}_{j\in\N}$
\mbox{to $\mu$.} Hence we may assume that $\partial
J\cap\partial\Sigma$ is a one--point set such that $\mu(\partial
J\setminus\partial\Sigma)=0$. Under these assumptions, there exist
monotone sequences of continuous functions $\{\varphi_K\}_{K\in\N}$
and $\{\psi_K\}_{K\in\N}$ such that $\varphi_K<\chi_{J}<\psi_K$ and
\mbox{$\int \psi_K-\varphi_K\,{\rm d}\mu<\frac{1}{K}$} for each
$K\in\N$. Since $\mu_j\overset{*}\to\mu$ (in the weak$^*$ topology)
as $j\to\infty$ and $\psi_K-\varphi_K$ is a continuous function, we
have that for each $K\in\N$ there exists $L_K\in\N$ such that
 \mbox{$\int \psi_K-\varphi_K\,{\rm d}\mu_j<\frac{2}{K}$} for all
 $j>
L_K$. It is easy to see that for each $K\in\N$ and for all $j>L_K$,
$$\left |\int \chi_J\,{\rm d}\mu-\int\chi_J\,{\rm d}\mu_j\right |<\frac{3}{K}+
\left |\int \varphi_K\,{\rm d}\mu-\int\varphi_K\,{\rm
d}\mu_j\right|. $$ This shows that
$\mu(J)=\lim_{j\to\infty}\mu_j(J)$.
\endproof

\end{lem}

\begin{lem}\label{nat} Let $\{x_{n_j}\}_{j=0}^{\infty}$ be a sequence in $\Sigma$ such that  $n_j\ge 1$ and
\mbox{$x_{n_j}\in{\rm dom}\,(P^{n_j-1})$} for all $j\in\N$. Any
accumulation point of the sequence of Borel probability measures
\begin{equation}\label{spm}
\mu_j=\frac{1}{n_j}\sum_{k=0}^{n_j-1}\delta_{P^k}(x_{n_j}),
\end{equation}
where $\delta_x$ is the Dirac probability measure on $\Sigma$
concentrated at $x$, is a non--atomic measure.
\end{lem}
\proof Let $\mu\in\mathscr{M}(\Sigma)$ be an accumulation point of
$\{\mu_j\}_{j\in\N}$. By taking a subsequence if necessary and by
renaming variables, we may assume that $\mu_j\overset{*}\to\mu$ as
$j\to\infty$. Since the set
\mbox{$D=\{z\in\Sigma\mid\mu(\{z\})>0\}$} is at most countable, for
each $p\in\Sigma$, there exists an open subinterval $I_p$ of
$\Sigma$ containing $p$ of length $\ell(I_p)$ arbitrarily small such
that $\mu(\partial I_p\setminus\partial\Sigma)=0$. By $(c)$ of
Theorem \ref{gst} and by Lemma 3.1 of Camelier--Gutierrez
\cite{CaGu}, for each $\epsilon>0$, there exist $\delta>0$ and
$N\in\N$, such that if $\ell(I_p)<\delta$ then for each $n\ge N$ and
$x\in{\rm dom}\,(P^{n-1})$,
$$\frac1n\sum_{k=0}^{n-1}\chi_{I_p}(P^k(x))<\epsilon.$$
Hence, for each $\epsilon>0$, there exist $\delta>0$ and $N\in\N$
such that if $\ell(I_p)<\delta$ then for all $j\ge N$,
$$\mu_j(I_p)=\frac{1}{n_j}\sum_{k=0}^{n_j-1}\chi_{I_p}(P^k(x_{n_j}))<\epsilon.$$
By Lemma \ref{Par}, for each $\epsilon>0$ there exists $\delta>0$
such that if $\ell(I_p)<\delta$ then
$$\mu(I_p)=\lim_{j\to\infty}\mu_j(I_p)\le\epsilon.$$
Hence, $\mu(\{p\})=0$ and so $\mu$ is non--atomic, which finishes
the proof.
\endproof

\begin{pro}\label{cphi} Suppose that there exist a constant $c>0$ and a continuous
function $\phi:\Sigma\to\R$ such that $\int \phi\,{\rm d}\nu<-c$ for
all $\nu\in\mathcal{M}_P(\Sigma)$. Then there exists $N\in\N$ such
that for each $n>N$ and for all $x\in{\rm dom}\,(P^{n-1})$,
$$\frac1n\sum_{k=0}^{n-1}\phi(P^k(x))<-c.$$
\end{pro}
\proof Assume by contradiction that there exists a sequence
$\{x_{n_j}\}_{j=0}^\infty\subset\Sigma$ such that for each $j\in\N$,
\mbox{$x_{n_j}\in{\rm dom}\,(P^{n_j-1})$} and
$$\frac{1}{n_j}\sum_{k=0}^{n_j-1}\phi(P^k(x_{n_j}))\ge -c.$$
The set $\mathscr{M}(\Sigma)$, endowed with the weak$^*$ topology, is a compact
metric space. Consequently, the sequence of Borel probability measures
$$\mu_{j}=\frac{1}{n_j} \sum_{k=0}^{n_j-1}\delta_{P^k(x_{n_j})},$$
has a subsequence that weakly$\sp{*}$ converges to a Borel
probability measure $\mu\in\mathscr{M}(\Sigma)$. By renaming
variables, we may assume that $\mu_j\overset{*}\to\mu$ as $j\to\infty$. By Lemma
\ref{Par} and by \mbox{Lemma \ref{nat},} we have that
$\mu(B)=\lim_{j\to\infty}\mu_j(B)$ for all Borel set
$B\in\mathcal{B}$. This combined with the fact that
$\lim_{j\to\infty}\mu_j(P^{-1}(B))=\lim_{j\to\infty}\mu_j(B)$ for
all Borel set $B$ yields that $\mu$ is $P-$invariant and so
$\mu\in\mathscr{M}_P(\Sigma)$. Since
the function $\phi$ is continuous, we have
$$\int\phi\,{\rm d}\mu=\lim_{j\to\infty}\int\phi\,{\rm d}\mu_j=\lim_{j\to\infty}\frac{1}{n_j}\sum_{k=0}^{n_j-1}\phi(P^k(x_{n_j}))
\ge -c,$$
by the definition of $\mu$ and by the way we have chosen the sequence $\{n_j\}_{j=0}^\infty$,
which contradicts the initial assumption that $\int
\phi\,{\rm d}\nu<-c$ for all $\nu\in\mathcal{M}_P(\Sigma)$. \endproof

\noindent{\bf Theorem B. }{\it Suppose that $X$ has divergence less
or equal to zero at its saddle--points and that $X$ has negative
Lyapunov exponents at a quasiminimal set $N$. Then $X$ has the
infinitesimal contraction property at $N$.} \proof Let $\Sigma_1$ be
a transversal segment to $X$ passing through a non--trivial
recurrent point $p\in N$ so small that the forward Poincar\'e Map
$T:\Sigma_1\to T(\Sigma_1)\subset C$ is a diffeomorphism, where $C$
is a special transversal circle for $N$. By the hypothesis on the
divergence of $X$ at its saddle--points every forward Poincar\'e Map
$P:\Sigma\to\Sigma$ induced by $X$ on a transversal segment $\Sigma$
to $X$ has properties $(P1)$ and $(P2)$ (see \cite{MSMM}). By the
hypothesis of negative Lyapunov exponents and by Lemma \ref{abc},
there exist a subsegment $\Sigma$ of $\Sigma_1$ passing through $p$
and a constant $c>0$ such that the forward Poincar\'e Map
$P:\Sigma\to\Sigma$ induced by $X$ satisfies $\int\log\vert
DP\vert\,{\rm d}\nu<-c$ for all $\nu\in\mathscr{M}_P(\Sigma)$. By
Lemma \ref{cext}, there exists a continuous function
$\phi:\Sigma\to\R$ everywhere defined, with $\log\vert
DP(x)\vert<\phi(x)$ for all $x\in{\rm dom}\,(P)\setminus
P^{-1}(\partial\Sigma)$, such that $\int\phi\,{\rm d}\nu<-c$ for all
$\nu\in\mathscr{M}_P(\Sigma)$. By Proposition \ref{cphi}, there
exists $N\in\N$ such that for all $n\ge N$ and for all $x\in{\rm
dom}\,(P^n)\setminus\mathcal{O}_n^-(x)$,
$$\frac1n\log\vert DP^n(x)\vert=\frac{1}{n}\sum_{k=0}^{n-1}\log\vert DP(P^k(x))\vert<\frac{1}{n}\sum_{k=0}^{n-1}\phi(P^k(x))<-c.$$
Thus $P^n$ is an infinitesimal contraction. By
 Proposition \ref{pro3}, $X$ has the infinitesimal \mbox{contraction} property at
 $N$.\bbox\vspace{0.2cm}

 To finish the paper, we now provide a sketch of the proof of
 Theorem C.\vspace{0.2cm}

 \noindent{\bf Theorem C.\,} {\it  Suppose that
$X\in\mathfrak{X}_H^r(M)$, $r\ge 2$, has the contraction property at
a {quasimi\-nimal} set $N$.  There exists $Y\in\mathfrak{X}^r(M)$
arbitrarily $C^r-$close to $X$ having one more saddle--connection
than $X$.}\vspace{0.2cm}

\noindent{\it Sketch of the proof.} In the smooth case, we may use
the same proof of Theorem \ref{saddle-connection} without any
changes. In the case in which $X\in\mathfrak{X}^r(M)$ we cannot use
that proof because in taking a $C^r-$flow box to make the
perturbation, the vector field so obtained is of class $C^{r-1}$.
Thus we have to make the perturbation directly on the surface (using
bump functions defined on the surface and using also the orthogonal
vector field to $X$) and to use the flow box coordinates only for
estimation purposes. \bbox
\medskip

\noindent{\bf Acknowledgments.} The authors are grateful to the professors A. L\'opez, B. Sc\'ardua, \mbox{D. Smania}, and
  M. A. Teixeira, who have made many useful remarks to a
previous version of this paper. We also would like to thank the professor Simon Lloyd
for checking the English. The first author was supported in
part by FAPESP Grant Tem\'atico \mbox{$\#$ 03/03107-9,} and by CNPq
Grant \mbox{$\#$ 306992/2003-5.} The second author was fully
\mbox{supported} by FAPESP Grants \mbox{$\#$ 03/03622-0} and $\#$
06/52650-5.

\bibliographystyle{alpha}

\end{document}